\documentclass[12pt]{amsart}
\usepackage{amsmath, amsthm, amscd, amsfonts}

\newtheorem{theorem}{Theorem}[section]

\theoremstyle{definition}
\newtheorem{definition}[theorem]{Definition}
\newtheorem{example}[theorem]{Example}
\newtheorem{problem}[theorem]{Problem}

\theoremstyle{remark}
\newtheorem{remark}[theorem]{Remark}

\numberwithin{equation}{section}
\begin{document}
\title{What is a Hilbert $C^*$-module?}
\author{Mohammad Sal Moslehian}
\address{Department of Mathematics, Ferdowsi University, P. O. Box 1159, Mashhad 91775, Iran.}
\email{moslehian@ferdowsi.um.ac.ir and moslehian@member.ams.org}


\begin{abstract} Hilbert $C^*$-modules are useful tools in $AW^*$-algebra theory, theory of operator algbras, operator K-theory, group representation theory and theory of operator spaces. The theory of Hilbert $C^*$-modules is very interesting on its own.
In this paper we give fundamentals of the theory of Hilbert
$C^*$-modules and examine some ways in which Hilbert $C^*$-modules
differ from Hilbert spaces.
\end{abstract}
\maketitle

\section{Introduction}

The notion of a Hilbert $C^*$-module is a generalization of the
notion of a Hilbert space. The first use of such objects was made by
I. Kaplansky \cite{KAP}. The research on Hilbert $C^*$-modules began
in the 70's in the work of the induced representations of
$C^*$-algebras by M. A. Rieffel \cite{RIE} and the doctoral
dissertation of W. L. Paschke \cite{PAS}. It is also used to study
Morita equivalence of $C^*$-algebras, KK-theory of $C^*$-algebras,
operator K-theory, $C^*$-algebra quantum group and theory of
operator spaces.In this paper we shall view fundamentals of the
theory of Hilbert $C^*$-module and discuss on some its aspects.

\section{Hilbert $C^*$-module}

\begin{definition} A pre-Hilbert module over a $C^*$-algebra ${\mathcal A}$ is a comlex
linear space ${\mathcal X}$ which is an algebraic right ${\mathcal
A}$-module equipped with an ${\mathcal A}$-valued inner product
$<.,.>: {\mathcal X}\times {\mathcal X}\to {\mathcal A}$ satisfying

(i) $<x,x>\ge 0, \textrm{~and~}<x,x>=0 \textrm{~iff~} x=0$;

ii) $<x,y+\lambda z> = <x,y>+\lambda <x,z>$;

(iii) $<x,ya>=<x,y>a$;

(iv) $<y,x>=<x,y>^*$

\noindent for each $x, y, z \in {\mathcal X}, \lambda \in
\mathbb{C}, a \in {\mathcal A}.$
\end{definition}

\begin{remark} Left pre-Hilbert $C^*$-modules are defined
similarly.
\end{remark}

\begin{example}

\noindent (i) Every inner product space is a left pre-Hilbert module
over $\mathbb{C}$.

\noindent (ii) If $I$ is a (closed) right ideal of ${\mathcal A}$,
then $I$ is a pre-Hilbert ${\mathcal A}$-module if we define
$<a,b>=a^*b$; in particular ${\mathcal A}$ is a pre-Hilbert
${\mathcal A}$-module.

\noindent (iii) Let $\{{\mathcal X}_i\}_{1\le i\le m}$ be a finite
family of pre-Hilbert ${\mathcal A}$-modules. Then the vector space
direct sum $\displaystyle{\bigoplus _{i=1} ^m {\mathcal X}_i}$ is a
pre-Hilbert ${\mathcal A}$-module if we
define$$(x_1,...,x_m)a=(x_1a,...,x_ma),
<(x_1,...,x_m),(y_1,...y_m)>=\displaystyle{\sum_{i=1}^m }<x_i,y_i>$$

\noindent (iv) If ${\mathcal H}$ is a Hilbert space, then algebraic
(vector space) tensor product ${\mathcal H} \otimes _{alg}{\mathcal
A}$ is a pre-Hilbert ${\mathcal A}$-module under$$(\xi\otimes
a).b=\xi\otimes (ab),~ <\xi\otimes a,\eta\otimes b>=<\xi,\eta>a^*b$$
Let $t=<\displaystyle{\sum_{i=1}^n}\xi_i\otimes
a_i,\displaystyle{\sum_{i=1}^n}\xi_i\otimes a_i>$ and let
$\{\epsilon_k\}_{1\le k\le m}$ be an orthogonal basis for the finite
dimensional subspace of ${\mathcal H}$ spanned by the $\xi_i$. If
$\xi_i=\displaystyle{\sum_{k=1}^m}\lambda_{ik}\epsilon_k$ , then
$$t=\sum_{k=1}^m(\sum_{i=1}^n\lambda_{ik}a_i)^*(\sum_{i=1}^n\lambda_{ik}a_i)\ge 0.$$
If $t=0$, then $\displaystyle{\sum_{i=1}^n}\lambda_{ik}a_i=0
\textrm{ ~for ~all~} k$, from which
$$\displaystyle{\sum_{i=1}^n\xi_i\otimes
a_i=\sum_i(\sum_k\lambda_{ik}\epsilon_k)\otimes
a_i=\sum_i\sum_k\epsilon_k\otimes\lambda_{ik}a_i=\sum_k\epsilon_k\otimes(\sum_i\lambda_{ik}a_i)=0}$$
\end{example}

\begin{theorem}If ${\mathcal X}$ is a Hilbert ${\mathcal A}$-module and $x,y\in
{\mathcal X}$, then $$<y,x><x,y>\le\|<x,x>\|<y,y>.$$
\end{theorem}
\begin{proof}
Suppose $\|<x,x>\|=1$. For $a\in {\mathcal A}$,we have
$0\le<xa-y,xa-y>=a^*<x,x>a-<y,x>a-a^*<x,y>+<y,y>\le
a^*a-<y,x>a-a^*<x,y>+<y,y>$ (since if c is a positive element of
${\mathcal A}$, then $a^*ca \le \| c \| a^*a \textrm{ ~for~ all~}
a\in {\mathcal A}$; \cite[Corollary I.4.6]{DAV}) Now put $a=<x,y>$
to get $a^*a\le<y,y>$.
\end{proof}

For $x\in {\mathcal X}$, Put $$\| x \|=\| <x,x> \|^{\frac{1}{2}}$$
Then $\|<x,y>\|^2=\|<x,y>^*<x,y>\|=\|<y,x><x,y>\| \le
\|<x,x>\|\|<y,y>\| =\| x\|^2\| y\|^2$. Hence
$$\|<x,y>\|\le\| x\|\| y\|$$
It follows that $<.,.>:{\mathcal X}\times {\mathcal X} \rightarrow
{\mathcal A}$ is continuous. Obviously $\| .\|$ is a norm on
${\mathcal X}$.

\begin{definition} A pre-Hilbert $C^*$-module ${\mathcal X}$ is called a
Hilbert $C^*$-module if ${\mathcal X}$ equipped with the above norm
is complete. Note that a Banach space completion of a pre-Hilbert
${\mathcal A}$-module, is a Hilbert $C^*$-module.
\end{definition}

We have $0\le<xa,xa>=a^*<x,x>a\le a^*\| x\|^2a$. Taking norm we get
$\|<xa,xa>\|\le\| a^*\| x\|^2a\|$ $\| xa\|^2\le\| x\|^2\|
a\|^2\textrm{ or}\| xa\|\le\| x\|\| a\|$. Hence ${\mathcal X}$ is a
normed ${\mathcal A}$-module. There is also an ${\mathcal A}$-valued
norm $\mid .\mid$ on ${\mathcal X}$ given by $\mid
x\mid=<x,x>^{1/2}$.

\begin{example} Let $\{{\mathcal X}_n\}$ be a sequence of Hilbert
${\mathcal A}$-modules. Then $\oplus_{i=1}^\infty {\mathcal
E}_n=\{\{x_n\}\mid x_n\in {\mathcal X}_n
,\sum_{n=1}^\infty<x_n,x_n>\textrm{~converges~in}~{\mathcal A}\}$
with operations $\{x_n\}+\lambda\{y_n\}=\{x_n+\lambda
y_n\}~,~\{x_n\}a=\{x_na\}~,
<\{x_n\},\{y_n\}>=\displaystyle{\sum_{n=1}^\infty}<x_n,y_n>$ is a
Hilbert $C^*$-module. Note that in $\displaystyle \bigoplus_{n=p}^q
{\mathcal X}_n$ we have
$$\|\sum_{n=p}^q<x_n,y_n>\|\le\|\sum_{n=p}^q<x_n,x_n>\|\|\sum_{n=p}^q<y_n,y_n>\|.$$
Hence $\displaystyle{\sum_{n=1}^\infty}<x_n,y_n>$ converges, by the
Cauchy criterion. For completeness, Let $\{u_k\}_k$ be a Cauchy
sequence in $\displaystyle{\bigoplus_{n=1}^\infty}{\mathcal X}_n$
and for all $k, ~ u_k=\{x_{n,k}\}_n$. Applying $\|
x_{n,k}-x_{n,l}\|=\|<x_{n,k}-x_{n,l},x_{n,k}-x_{n,l}>\|\le
\|\displaystyle{\sum_{n=1}^\infty}<x_{n,k}-x_{n,l},x_{n,k}-x_{n,l}>\|=\|
u_k-u_l\|^2$, we conclude that $\{x_{n,k}\}_k$ is Cauchy, for each
$n\in N$. Hence for each $n$, there exists $v_n$ such that
$\lim_kx_{n,k}=v_n$. Now put $v=\{v_n\}.~\textrm{ Then~}
\lim_ku_k=v.$
\end{example}

Hilbert $C^*$-modules behave like Hilbert spaces in some way, for example:
$$\| x\|=\sup \{\|<x,y>\|~;y\in {\mathcal X},\| y\|\le 1\}$$
But there is one fundamental way in which Hilbert $C^*$-modules
differ from Hilbert spaces: Given a closed submodule ${\mathcal F}$
of a Hilbert ${\mathcal A}$-module ${\mathcal X}$, define
$${\mathcal F}^{\bot}=\{y\in {\mathcal X}\mid<x,y>=0, \forall x\in {\mathcal F}\}.$$
Then ${\mathcal F}^{\bot}$ is a closed submodule. But usually
${\mathcal X}\ne {\mathcal F}+{\mathcal F}^{\bot}~,~{\mathcal
F}^{\bot\bot}\ne {\mathcal F}.$

\begin{example} Take ${\mathcal A}=C([0,1])$, ${\mathcal X}={\mathcal A}$ and ${\mathcal F}=\{f \in {\mathcal A}\mid
f(1/2)=0\}$. Then~${\mathcal F}^{\bot}=\{ 0 \}$. $(\forall g\in
{\mathcal F}^{\bot};g(t)\mid t-1/2\mid=0.~\textrm{Hence}~g\equiv0$
by continuity) Pythagoras' equality stating $\xi,\eta\in {\mathcal
H}$ and $\xi\bot\eta$ imply $\|\xi+\eta\|^2=\|\xi\|^2+\|\eta\|^2$,
does not hold, in general, for Hilbert $C^*$-modules. For example
consider ${\mathcal A}=C([0,1]\cup[2,3])$ as a Hilbert ${\mathcal
A}$-module,
\[f(x)= \left \{
\begin{array}{cc}1& x \in [0,1]\\0& x \in [2,3] \end{array} \right .
\textrm{\,\,\,\,and\,\,\,\,} g(x)= \left \{ \begin{array}{cc}0& x
\in [0,1]\\1& x \in [2,3] \end{array} \right .\] Then
$$<f,g>=fg=0,~ \| f+g\|=1,\textrm{ ~and~} \| f\|=\| g\|=1$$
\end{example}

\begin{definition} A closed submodule ${\mathcal F}$ of a Hilbert ${\mathcal A}$-module
${\mathcal X}$ is topologically complementable if there is a closed
submodule ${\mathcal G}$ such that ${\mathcal F}+{\mathcal
G}={\mathcal X}$ and ${\mathcal F}\cap {\mathcal G}=\{0\}$.
${\mathcal F}$ is called orthogonally complemented if we have the
further condition ${\mathcal F}\bot {\mathcal G}$.
\end{definition}

Not every topologically complemented is orthogonally complemented as
the following example shows.

\begin{example} Let ${\mathcal A}=C([0,1]),~J=\{f\in {\mathcal A}\mid f(0)=0\}~\simeq
C_0((0,1])$ and ${\mathcal X}={\mathcal A} \oplus J$ as a Hilbert
${\mathcal A}$-module. If ${\mathcal F}=\{(f,f)\mid f\in J\}$, then
${\mathcal F}^\bot =\{(g,-g)\mid g\in J\}~, {\mathcal F}+{\mathcal
F}^\bot =J+J\ne {\mathcal X}$ and ${\mathcal G}=\{(f,0)\mid f\in
{\mathcal A} \}$ is a topological complement for ${\mathcal F}$.
\end{example}

\section{Full Hilbert $C^*$-modules}

Let ${\mathcal X}$ be a Hilbert ${\mathcal A}$-module and let
$\{e_{\lambda}\}$ be an approximate unit for ${\mathcal A}$. For
$x\in {\mathcal
E}~,~<x-xe_{\lambda},x-xe_{\lambda}>=<x,x>-e_{\lambda}<x,x>-<x,x>e_{\lambda}+e_{\lambda}<x,x>e_{\lambda}\displaystyle
\stackrel{\lambda}{\to}0$ Hence $\lim_\lambda xe_\lambda=x$. It
follows that ${\mathcal X}{\mathcal A}$ which is defined to be the
linear span of $\{xa\mid x\in {\mathcal X}, a\in {\mathcal A} \}$ is
dense in ${\mathcal X}$ and if ${\mathcal A}$ is unital, then
$x.1=x$. Clearly $<{\mathcal X},{\mathcal X}>= \textrm{ span}
\{<x,y>\mid x,y\in {\mathcal X}\}$ is a $*$-bi-ideal of ${\mathcal
A}$.

\begin{definition} If $<{\mathcal X},{\mathcal X}>$ is dense in ${\mathcal A}$, then ${\mathcal X}$ is called
full. For example ${\mathcal A}$ as an ${\mathcal A}$-module is
full.
\end{definition}

\section{Adjoinable Maps between Hilbert $C^*$-modules}

Suppose that ${\mathcal X}$ and ${\mathcal F}$ are Hilbert
${\mathcal A}$-modules and
$$L({\mathcal X},{\mathcal F})=\{t:{\mathcal X}\to {\mathcal F}\mid\exists t^*:{\mathcal F}\to
{\mathcal X};<tx,y>=<x,t^*y>\}$$ $t$ must be ${\mathcal A}$-linear,
since $<t(xa),y>=<xa,t^*y>=a^*<x,t^*y>=a^*<tx,y>=<(tx)a,y>$ for all
$y$. So $<t(xa)-(tx)a ,y>=0$ and hence
$<t(xa)-(tx)a,t(xa)-(tx)a>=0$. It follows that $t(xa)-(tx)a=0$
Similarly; $t(\lambda x+y)=\lambda tx+ty$ t must be bounded since
For each $x$ in the unit ball of ${\mathcal X}$, define
$f_x:{\mathcal F} \to {\mathcal A}$ by $f_x(y)=<tx,y>=<x,t^*y>$.
Then $\| f_x(y)\|\le\| x\| \| t^*y\|\le\| t^*y\|$. So $\{\| f_x\| ;
x\in {\mathcal X}_1\}$ is bounded. This and $\|
tx\|=\displaystyle{\sup_{y\in {\mathcal F}_1}\|<tx,y>\|=\sup_{y\in
{\mathcal F}_1}}\| f_x(y)\|=\| f_x\|$ show that t is bounded.
$L({\mathcal X},{\mathcal F})$ is called the space of adjointable
maps and we put $L({\mathcal X})=L({\mathcal X},{\mathcal X})$.

A bounded ${\mathcal A}$-linear map need not be adjointable:

\begin{example} Let ${\mathcal F}={\mathcal A} =C([0,1]),~{\mathcal X}=\{f\in {\mathcal A} ; f(1/2)=0\}$ and
$i: {\mathcal X} \to {\mathcal F}, f \mapsto f$ be the inclusion
map. If $i$ were adjointable and $1$ denote the identity element of
${\mathcal A}$, then for all $x\in {\mathcal X};
<x,i^*(1)>=<i(x),1>=<x,1>$. Hence $i^*(1)=1.$ But $1\notin {\mathcal
E}$ and so $i$ cannot be adjointable. (Collate this with Hilbert
spaces.)
\end{example}

\section{$L({\mathcal X})$ as a $C^*$-algebra}

If $t\in L({\mathcal X},{\mathcal F})$, then $t^*\in L({\mathcal
F},{\mathcal X})$. If ${\mathcal G}$ is a Hilbert ${\mathcal
A}$-module and $s\in L({\mathcal F},{\mathcal G})$, then $st\in
L({\mathcal X},{\mathcal G})$. Hence $L({\mathcal X})$ is a
$*$-algebra. If $t_n\to t$ then
\begin{eqnarray*}\|
t^*_ny-t^*_my \|&=&\sup_{x\in {\mathcal X}_1} \|<x,(t^*_n-t^*_m)y>
\|\\
&=&\sup_{x\in {\mathcal X}_1}\|<(t_n-t_m)x,y>\| \\
&\le& \sup_{x\in {\mathcal X}_1}\|(t_n-t_m)x\|\| y\|\\
&\le& \| t_n-t_m\|\| y\|.
\end{eqnarray*}
It follows that
$\{t^*_ny\}$ converges to, say, $sy$. Hence
$$<tx,y>=\lim_n<t_nx,y>=<x,\lim_nt^*_ny>=<x,sy>.$$
Thus $t\in L({\mathcal X})$. Thus $L({\mathcal X})$ is a closed
subset of $$B({\mathcal X})=\{T:{\mathcal X}\to {\mathcal X}\mid
T\textrm{ ~is~ linear ~and~ bounded}\}.$$ Hence $L({\mathcal X})$ is
a Banach algebra. Moreover,
\begin{eqnarray*}
\|t\|^2&=&\sup_{x\in {\mathcal X}_1}\|tx\|^2\\
&=& \sup_{x\in {\mathcal X}_1}\|<tx,tx>\|\\
&=& \sup_{x\in {\mathcal X}_1}\|<t^*tx,x>\|\\
&\le& \| t^*t\|\\
&\le& \| t^* \| \| t\|\\
&\le& \| t\|^2.
\end{eqnarray*} Hence $\|
t\|^2=\| t^*t\|$. Thus $L({\mathcal X})$ is a $C^*$-algebra.

\section{``Compact'' operators}
Let ${\mathcal X}$ and ${\mathcal F}$ be Hilbert ${\mathcal
A}$-modules. For $x\in {\mathcal X}$ and $y\in {\mathcal F}$ define
$\Theta_{x,y}: {\mathcal F} \to {\mathcal X}$ by
$\Theta_{x,y}(z)=x<y,z>$. Then $\Theta^*_{x,y}=\Theta_{y,x}
\Theta_{x,y}\Theta_{u,v}=\Theta_{x<y,u>,v},
t\Theta_{x,y}=\Theta_{tx,y}$. Let $K({\mathcal F},{\mathcal X})$,
the set of ``compact'' operators, be the closed linear span of
$\{\Theta_{x,y}\mid x\in {\mathcal X},y\in {\mathcal F}\}$. It is
straightforward to show that $K({\mathcal X})=K({\mathcal
E},{\mathcal X})$ is a closed bi-ideal of $L({\mathcal X})$.

\begin{example} In the case ${\mathcal X}={\mathcal A}$, we have $K({\mathcal X})\simeq {\mathcal A}$ given
by $\Theta_{a,b}\longleftrightarrow ab^*$; cf. \cite[Page 10]{LAN}.
If ${\mathcal A}$ is unital, then $K({\mathcal A})=L({\mathcal A})$.
In fact if $t\in L({\mathcal A})$, then
$tx=t(1.x)=t(1)x=\Theta_{t(1),1}(x)$ and so $t=\Theta_{t(1),1}\in
K({\mathcal A})$. If ${\mathcal A}$ is not unital, $K({\mathcal
A})\simeq {\mathcal A}$ and $L({\mathcal A})\simeq M({\mathcal A})$
, the multiplier of ${\mathcal A}$ , the maximal $C^*$-algebra
containing ${\mathcal A}$ as an essential ideal.
\end{example}

\begin{theorem}
If ${\mathcal X}$ is a Hilbert ${\mathcal A}$-module, then
$M(K({\mathcal X}))\simeq L({\mathcal X})$; \textrm{cf.
\cite[Theorem 15.2.12]{WEG}}
\end{theorem}

\section{Dual of a Hilbert $C^*$-module}

For a pre-Hilbert ${\mathcal A}$-module ${\mathcal X}$, let
${\mathcal X}^{\#}$, the dual of ${\mathcal X}$, be the set of all
bounded ${\mathcal A}$-linear maps ${\mathcal X}$ into ${\mathcal
A}$. ${\mathcal X}^{\#}$ is a linear space and a right ${\mathcal
A}$-module under $(\alpha+\beta)(x)=\alpha(x)+\beta(x) ,
(\lambda\alpha)(x)=\-\lambda\alpha(x), (\alpha.a)(x)=a^*\alpha(x),
\quad (a \in{\mathcal A}, \lambda\in {\mathbb C}, x\in {\mathcal
E}$.

\noindent $\psi: {\mathcal X} \to {\mathcal X}^{\#}, x \mapsto
\widehat{x}$ , where $\widehat{x}: {\mathcal X} \to {\mathcal A},
\widehat{x}(y)=<x,y>$ is an isometric ${\mathcal A}$-linear map.
Identify ${\mathcal X}$ with $\widehat{{\mathcal X}}=\{\widehat{x};
x\in {\mathcal X}\}$ as an submodule of ${\mathcal X}^{\#}$.

\begin{problem} Could $<.,.>$ on ${\mathcal X}$ be extended to an ${\mathcal A}$-valued
inner product on ${\mathcal X}^{\#}$? W. L. Paschke showed that it
can be done if ${\mathcal A}$ is a von Neumann algebra \cite{PAS}.
\end{problem}

\section{Self-dual Hilbert $C^*$-modules}

\begin{definition} ${\mathcal X}$ is called self dual if
$\widehat{{\mathcal X}}={\mathcal X}^{\#}$.
\end{definition}

\begin{theorem} ${\mathcal A}$ as a Hilbert ${\mathcal A}$-module is self-dual iff ${\mathcal A}$ is unital.
\end{theorem}
\begin{proof} Let ${\mathcal A}$ be unital with unit $1$ and let $t\in {\mathcal A}^{\#}$. Then for all $a\in
{\mathcal A}$;
$$t(a)=t(1.a)=t(1).a=<t(1)^*,a>=(t(1)^*)^\wedge (a).$$ Hence
$t=(t(1)^*)^\wedge \in \widehat{{\mathcal A}}\subseteq {\mathcal
A}^{\#}$. If $\widehat{{\mathcal A}}={\mathcal A}^{\#}$ , then $i:
{\mathcal A} \to {\mathcal A}, i(x)=x$ being bounded ${\mathcal
A}$-linear is of the form $\widehat{a}$ for some $a\in {\mathcal
A}$. Hence for all $x\in {\mathcal A},
x=i(x)=\widehat{a}(x)=<a,x>=a^*x$. So $a^*$ is the unit of
${\mathcal A}.\Box$
\end{proof}

\begin{theorem} If ${\mathcal F}$ is a self-dual closed ${\mathcal A}$-submodule
of a Hilbert $C^*$-algebra ${\mathcal X}$, then ${\mathcal F}$ is
orthogonally complementable.
\end{theorem}
\begin{proof} For each $e\in {\mathcal X}$ there exists some $f\in {\mathcal F}$ such that $<e,.>=<f,.>$
on ${\mathcal F}$. So $<e-f,.>=0$ on ${\mathcal F}$. Hence $e-f\in
{\mathcal F}^\bot$. Since $e=(e-f)+f$, ${\mathcal F}$ is
orthogonally complementable.
\end{proof}

Saworotnow showed that every Hilbert ${\mathcal A}$-module over a
finite dimensional $C^*$-algebra ${\mathcal A}$ is self-dual
\cite{SAW}. If ${\mathcal X}$ is a Hilbert ${\mathcal A}$-module and
${\mathcal A}$ is unital, definie $\widetilde{{\mathcal X}}$ to be
the set of all bounded ${\mathcal A}$-module maps from ${\mathcal
E}$ to ${\mathcal A}^{**}$. Clearly $\widetilde{{\mathcal X}}$ is
right ${\mathcal A}^{**}$-module. and the map $x\mapsto\widehat{x}$
is an embedding of ${\mathcal X}$ in $\widetilde{{\mathcal X}}$. As
Paschke demonstrate, there is an extension of the inner product of
${\mathcal X}$ to an ${\mathcal A}^{**}$-valued inner product of
$\widetilde{{\mathcal X}}$ for which $\widetilde{{\mathcal X}}^{\#}
\simeq \widetilde{{\mathcal X}}$; \cite{PAS}. Thus every Hilbert
module can be embedded in a self-dual Hilbert module.

\section{$C^*$-reflexivity}

Denote by $ {\mathcal X}^{\#\#}$ the ${\mathcal A}$- module of all
bounded ${\mathcal A}$- module maps from ${\mathcal X}^{\#}$ into
${\mathcal A}$. Let $\Omega$ be the module map $\Omega: {\mathcal
E}\to {\mathcal X}^{\#\#},~ \Omega(x)(\tau)=\tau(x)^*; x\in
{\mathcal X},\tau\in {\mathcal X}^{\#}$.

\begin{definition} A Hilbert ${\mathcal A}$- module ${\mathcal X}$ is called
$C^*$-reflexive , if $\Omega$ is a module isomorphism. In the
Hilbert $W^*$-module setting, $C^*$-reflexivity is equivalent to the
self-duality; \cite{PAS}
\end{definition}

\section{Riesz Theorem for Hilbert $C^*$- modules }

Let ${\mathcal X}$ be a Hilbert ${\mathcal A}$-module and $x \in
{\mathcal X}$.
$$\Theta_{a,z}(y)=a<z,y>=<za^*,y>=(za^*)^\wedge(y)(*)$$
and so $\widehat{x}\in K({\mathcal X},{\mathcal A})$, whenever $x$
is of the form $za^*$. Since ${\mathcal X}{\mathcal A}$ is dense in
${\mathcal X}$ for each $x\in {\mathcal X}$ there is a sequence
$\{x_n\}$ in ${\mathcal X}{\mathcal A}$ such that
$\displaystyle{\lim_nx_n=x}$ . But ${\mathcal X} \to L({\mathcal
E},{\mathcal A}), x \mapsto \widehat{x}$ is comtinuous (isometry).
Hence $\widehat{x}=\displaystyle{\lim_n\widehat{x}_n} \in
K({\mathcal X},{\mathcal A})$, since $K({\mathcal X},{\mathcal A})$
is closed in $L({\mathcal X},{\mathcal A})$.

(*) shows that every element of $K({\mathcal X},{\mathcal A})$ is of
the from $\widehat{x}=<x,.>$ for some $x\in {\mathcal X}$, a Riesz
theorem for Hilbert $C^*$-modules.

\section{Linking Algebra}

Let $\Theta_{x,y}$ be denoted by $[x,y]$. Suppose ${\mathcal X}$ is
a Hilbert ${\mathcal A}$-module. Then $[.,.]$ is an $L({\mathcal
E})$-valued inner product on ${\mathcal X}$, for which ${\mathcal
E}$ becomes a left Hilbert $L({\mathcal X})$-module. A Hilbert
${\mathcal A}$-module ${\mathcal X}$ can be embedded into a certain
$C^*$-algebra $\Lambda({\mathcal X})$:

Denote by ${\mathcal F}={\mathcal X} \oplus {\mathcal A}$, the
direct sum of Hilbert ${\mathcal A}$-modules ${\mathcal X}$ and
${\mathcal A}$. Recall $<(x_1,a_1),(x_2,a_2)>=<x_1,x_2>+a^*_1a_2$.
Identify each $x\in {\mathcal X}$ with ${\mathcal A} \to {\mathcal
E}, a \mapsto xa$. The adjoint of this map is $x^*(y)=<x,y>$,
(indeed $x^*=\widehat{x}$). Put
$$\Lambda({\mathcal X})= \{ \left[ \begin{array}{cc}b&x\\y^*&a
\end{array}\right] \mid a\in {\mathcal A}, x,y\in {\mathcal X}, b\in L({\mathcal X})\}.$$
$\Lambda({\mathcal X})$ is a $C^*$-subalgebra of $L({\mathcal F})$,
called the linking algebra of ${\mathcal X}$. Then ${\mathcal X}
\simeq \left[
\begin{array}{cc}0&{\mathcal X}\\0&0 \end{array}\right],~{\mathcal A} \simeq \left[
\begin{array}{cc}0&0\\0&{\mathcal A} \end{array}\right]$ and $L({\mathcal X}) \simeq
\left[ \begin{array}{cc}L({\mathcal X})&0\\0&0 \end{array}\right]$.
Furthermore,$<x,y>$ of ${\mathcal X}$ becomes the product $x^*y$ in
$\Lambda({\mathcal X})$ and the module multiplication ${\mathcal
E}\times {\mathcal A}\to {\mathcal X}$ becomes a part of the
internal multiplication of $\Lambda({\mathcal X})$. Suppose now that
$\pi: \Lambda({\mathcal X})\to B({\mathcal H})$ is a faithful
representation of $\Lambda({\mathcal X})$. Then by restriction,
$\pi$ gives maps $\phi=\pi\mid_{\mathcal A}$ and
$\Phi=\pi\mid_{\mathcal X}$ such that $\phi(<x,y>)=\Phi(x)^*\Phi(y)$
and $\Phi(xa)=\Phi(x)\phi(a), (x,y\in {\mathcal X}, a\in {\mathcal
A})$ Hence ${\mathcal X}$ can be isometrically embedded into
$B({\mathcal H})$ via $\Phi$. Therefore ${\mathcal X}$ can be
regarded as an operator space; i.e. a closed subspace of
$B({\mathcal H})$ for some Hilbert space ${\mathcal H}$.


\begin{thebibliography}{99}
\bibitem{DAV} K. R. Davidson, \textit{$C^*$-algebras by Example}, Fields Ins. Monog., 1996.
\bibitem{FRA1} M. Frank, \textit{Self-duality and $C^*$-reflexivity of Hilbert $C^*$-modules}, Zeitschr Anal. Anw. \textbf{9}(1990), 165--176.
\bibitem{FRA2} M. Frank, \textit{A multiplier approach to the Lance-Blecher theorem}, J. for Anal and its Appl., \textbf{16} (1997), no. 3, 565--573.
\bibitem{FRA3} M. Frank, \textit{Geometrical aspects of Hilbert $C^*$-modules}, Positivity, 1998.
\bibitem{KAP} I. Kaplansky, \textit{Modules over operator algebras}, Amer J. Math., \textbf{75} (1953), 839--858.
\bibitem{LAN} E. C. Lance, \textit{Hilbert $C^*$-modules}, LMS Lecture Note Series 210, Cambridge University Press, Cambridge, 1995.
\bibitem{MAG} B. Magajna, \textit{Hilbert modules and tensor products of operator spaces},  Linear operators (Warsaw, 1994), 227--246, Banach Center Publ., 38, Polish Acad. Sci., Warsaw, 1997.
\bibitem{PAS} W. L. Paschke, \textit{Inner product modules over B*-algebras}, Trans Amer. Math. Soc. \textbf{182} (1973), 443--468.
\bibitem{RIE} M. A. Rieffel, \textit{Morita equivalence representations of $C^*$-algebras}, Adv. in Math., \textbf{13} (1974), 176--257.
\bibitem{SAW} P. P. Saworotnow, \textit{A generalized Hilbert space}, Duke Math. J. \textbf{35} (1968), 191--197.
\bibitem{WEG} N. E. Wegge-Olsen, \textit{K-theory and $C^*$-algebras}, Oxford University Press, New York, 1993.
\end{thebibliography}
\end{document}